# Geometric Constructions of Regular Polygons and Applications to Trigonometry

José Gilvan de Oliveira, Moacir Rosado Filho, Domingos Sávio Valério Silva

*Abstract*: In this paper, constructions of regular pentagon and decagon, and the calculation of the main trigonometric ratios of the corresponding central angles are approached. In this way, for didactic purposes, it is intended to show the reader that it is possible to broaden the study of Trigonometry by addressing new applications and exercises, for example, with angles 18°, 36° and 72°. It is also considered constructions of other regular polygons and a relation to a construction of the regular icosahedron.

## Introduction

The main objective of this paper is to approach regular pentagon and decagon constructions, together with the calculation of the main trigonometric ratios of the corresponding central angles. In textbooks, examples and exercises involving trigonometric ratios are usually restricted to so-called notable arcs of angles 30°, 45° and 60°. Although these angles are not the main purpose of this article, they will be addressed because the ingredients used in these cases are exactly the same as those we will use to build the regular pentagon and decagon. The difference to the construction of these last two is restricted only to the number of steps required throughout the construction process. By doing so, we intend to show the reader that it is possible to broaden the study of Trigonometry by addressing new applications and exercises, for example, with angles 18°, 36° and 72°. The ingredients used in the text are those found in Plane Geometry and are very few, restricted to intersections of lines and circumferences, properties of triangles, and Pythagorean Theorem.

For each integer $n \geq 3$, let $\wp_n$ be a regular $n$-sided polygon inscribed on a circle of unit radius and center $O$, considered on an established plane. Let $l_n$ be the side length of $\wp_n$ and let $\theta_n$ be the inner center angle in $O$ of a triangle with vertices $O$ and two consecutive vertices of $\wp_n$. So, since the unit radius circumference is divided into $n$ arcs of the same measure, the center angle $\theta_n$ is given by $\theta_n = 360°/n$, and hence, as in Figure 1, we have $\sin(\theta_n/2) = l_n/2$ and, therefore, $\sin(180°/n) = l_n/2$. One question we will study is how to effectively divide the unit circumference in this way, i.e., into $n$ equal parts. The main objective of this paper is to provide an answer to cases $n = 5$ and $n = 10$, i.e., to construct the regular pentagon $\wp_5$ and the regular decagon $\wp_{10}$. We will do this by repeating the same techniques as the previous two cases $n = 3$ and $n = 4$. As a consequence directly related to the construction of the regular pentagon, we will also get the cases $n = 10$ and $n = 20$.

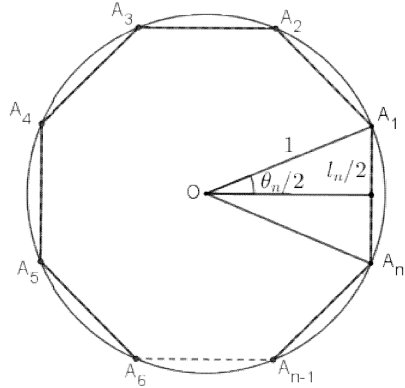

Figure 1

## Equilateral Triangle and Square

1) Case $n = 3$.

In this case, we have $\theta_3 = 360°/3 = 120°$ and $\sin(\theta_3/2) = \sin 60° = l_3/2$. Consider Figure 2 below. Given the line segment $AA'$, we determine its midpoint $O$, obtained from the intersection of the two circumferences with centers at the endpoints of segment $AA'$ and the same radius $r$ greater than or equal to half the length of $AA'$. The two triangles with common side $OA$ are isosceles and congruent, by the case $LLL$ of triangles congruence. Similarly, considering segment $OA$, we obtain its midpoint $C$. In addition, $C$ is also the midpoint of segment $BB'$. In particular, if $r = OA$, then the two triangles are equilateral, the measurement of angle $\alpha$ is $60°$ and the two triangles with common side $OC$ are congruent, by the case $LAL$. Therefore, they are also right triangles. Finally, since we have $\alpha$ as an exterior angle of triangle $A'OB$, the measurement of the base angles of this isosceles triangle is $\alpha/2 = 30°$. The same conclusion holds for the triangle $A'OB'$. This allows us to conclude that $A'B'B$ is an equilateral triangle and thus the angle $\widehat{OBC}$ is equal to $\alpha/2 = 30°$. Moreover, considering the segment $OA$ as unitary, by the Pythagorean Theorem applied to the right triangle OCB, we conclude that $\sin 30° = 1/2$, $\sin 60° = \sqrt{3}/2$ and that the side $l_3$ of the equilateral triangle $A'B'B$ is $l_3 = \sqrt{3}$.

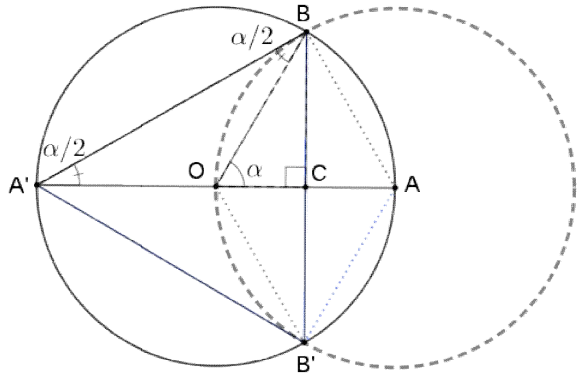
Figure 2

It is noteworthy that the procedure described, taking the intersection of two particular circumferences, effectively divides a circumference into three equal parts and thus constructs the regular polygon $\wp_3$. The process can be summarized in two steps: a) Given segment $AA'$ we determine its midpoint $O$. b) Then, from $O$, determine the intersection points $B$ and $B'$ of the segment $OA$ with the circumference of center $O$ and radius $OA$. Another highlight is that the method effectively allows obtaining the mediatrix of any segment and, therefore, bisector, median and height, relative to the base, of any isosceles triangle. This will be used in the other cases considered below. In the case of Figure 2 we have these concepts highlighted, for example, both for triangle $OAB$ and triangle $A'B'B$. Thus, point $C$ is midpoint of segment OA as well of segment $BB'$.

2) Case $n = 4$.

In this case, we have $\theta_4 = 360°/4 = 90°$ and $\sin(\theta_4/2) = \sin 45° = l_4/2$. Consider Figure 3 below.

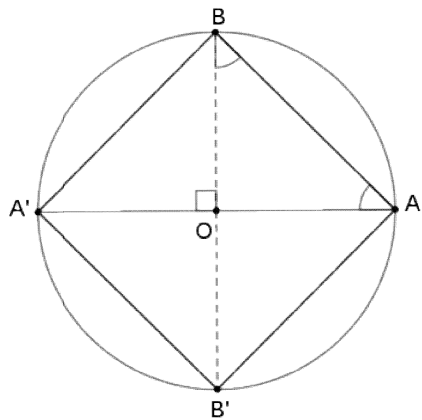
Figure 3

Given the line segment $AA'$, following the procedure described in the previous case, we determine its midpoint $O$. Then we determine the intersection points $B$ and $B'$ of the mediatrix of $AA'$ with the circumference of center $O$ and radius $OA$. Thus, we obtain the right triangle $AOB$. Since this triangle is also isosceles, considering $OA$ to be unitary and using the Pythagorean Theorem, we thus effectively obtain the inscribed square $\wp_4$. Also, $\sin 45° = \sqrt{2}/2$ and $l_4 = \sqrt{2}$.

## Regular Pentagon and Decagon

3) Case $n = 5$.

In this case, we have $\theta_5 = 360°/5 = 72°$ and $\sin(\theta_5/2) = \sin 36° = l_5/2$. Hence, from the sum of the inner angles of a triangle and the equalities $180° = 36° + 2 \cdot 72°$ e $180° = 36° + (36° + 36°) + 72°$, this suggests the existence of two similar isosceles triangles with bases 1 and 1, as in Figure 4 below, where $\alpha = 36°$. Such triangles will subsequently be obtained effectively, in which case, $l$ is the $l_{10}$ side of the regular polygon $\wp_{10}$. It follows from the similarity of these triangles, that $l = l/1 = (1-l)/l$, and therefore, $l = (\sqrt{5} - 1)/2$, which is the famous golden ratio ([1] and [4]).

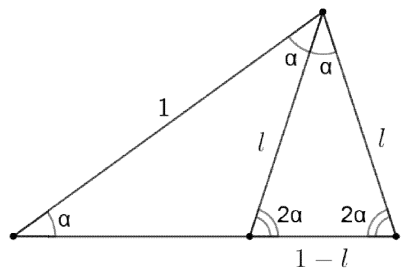

Figure 4

With knowledge of these facts, we will now effectively construct the regular polygon $\wp_5$, i.e., the regular pentagon. Consider Figure 5 below.

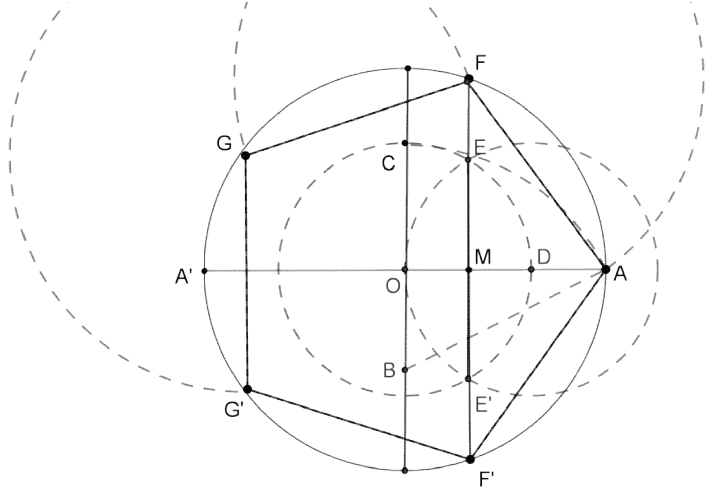

Figure 5

Given the line segment $AA'$ we determine, using the procedure detailed in the first case treated ($n = 3$), its mediator containing the midpoint $O$ of $AA'$ and also the circumference with center $O$ and radius $OA$. Then we determine point $B$ of the aforementioned mediatrix so that $OB$ is half the radius of the circumference. Assuming that the segment $OA$ is unitary, by the Pythagorean Theorem applied to the right triangle $AOB$, we have $AB = \sqrt{5}/2$. Therefore, $OC = AB - 1/2 = (\sqrt{5} - 1)/2$, which is the golden ratio, where $C$ is one of the intersection points of the same mediatrix with the circumference with center $B$ and radius $AB$. Point D is the intersection of segment $OA$ with the circumference with center $O$ and radius $OC$. Finally, points $F$ and $F'$ are the intersection points of the circumference of center $O$ and radius 1 with the mediatrix of segment $OD$. So the points $F'$, $A$ and $F$ are three consecutive points of $\wp_5$. The remaining points of $\wp_5$ are the points $G$ e $G'$ obtained from the intersection of the unitary circumference of center $O$ with the circumferences of the same radius $AF$ and centers at the points $F$ e $F'$. Also, as the base angles of the isosceles triangle $OFD$ are equal to 72° and the height $FM$ of this triangle relative to the vertex $F$ is half of $FF'$, it follows from Pythagorean Theorem that

$$\sin 72° = FM = \sqrt{(5 + \sqrt{5})/8} = \sqrt{(10 + 2\sqrt{5})}/4,$$

$$FF' = \sqrt{(10 + 2\sqrt{5})}/2,$$

and using the bisector of the triangle $OFD$ relative to vertex $F$, we have

$$\sin 18° = OM/OF = (\sqrt{5} - 1)/4.$$

Then, calculating the height of the triangle $OFD$ relative to vertex $D$ using the Pythagorean Theorem (see Figure 6), we obtain

$$\sin 36° = DH/DF = \sqrt{(5 - \sqrt{5})/8} = \sqrt{(10 - 2\sqrt{5})}/4.$$

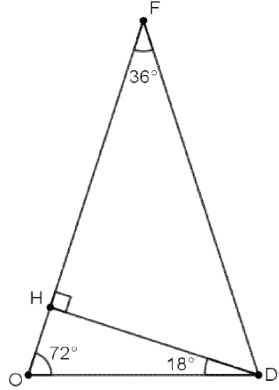

Figure 6

Hence, the sides of the regular pentagon $\wp_5$ and the regular decagon $\wp_{10}$ are, respectively

$$l_5 = \sqrt{(10 - 2\sqrt{5})/2} \text{ and } l_{10} = (\sqrt{5} - 1)/2.$$

Another interesting property about the regular pentagon, as can be seen from the above informations, is the golden ratio as the ratio of the side length to the diagonal of the regular pentagon. In particular, using the $LLL$ case of triangle similarity, if the $AF$ side in Figure 5 is 1, then the diagonal $FF'$ is $\frac{1+\sqrt{5}}{2}$ (see [4]).

From previous results and the relationship $\cos^2 \theta + \sin^2 \theta = 1$, we have the tables.

| $\theta$ | 30° | 45° | 60° |
|---|---|---|---|
| $\sin \theta$ | $\frac{1}{2}$ | $\frac{\sqrt{2}}{2}$ | $\frac{\sqrt{3}}{2}$ |
| $\cos \theta$ | $\frac{\sqrt{3}}{2}$ | $\frac{\sqrt{2}}{2}$ | $\frac{1}{2}$ |

| $\theta$ | 18° | 36° | 72° |
|---|---|---|---|
| $\sin \theta$ | $\frac{\sqrt{5}-1}{4}$ | $\frac{\sqrt{10-2\sqrt{5}}}{4}$ | $\frac{\sqrt{10+2\sqrt{5}}}{4}$ |
| $\cos \theta$ | $\frac{\sqrt{10+2\sqrt{5}}}{4}$ | $\frac{\sqrt{5}+1}{4}$ | $\frac{\sqrt{5}-1}{4}$ |

**Example**: The Master Plan of a certain coastal city, in order to ensure proper beach insolation at a certain date and time of the year, states that the height of the buildings must not exceed the angle of 36°, measured from the boardwalk along a horizontal line, to the highest point of the facade, as shown in Figure 7 below. Assuming that the distance from the boardwalk to the base of a building is 21 meters, what is the maximum height $h$ of the building?

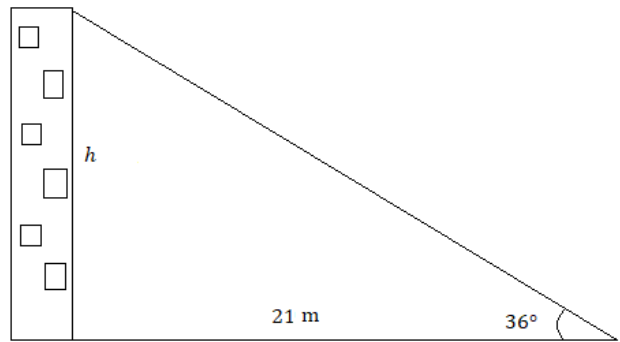

Solution: For the maximum height $h$ of the building, we have

$$\frac{h}{21} = \tan 36° = \frac{\sin 36°}{\cos 36°} = \frac{\sqrt{10 - 2\sqrt{5}}}{1 + \sqrt{5}}.$$

Hence, $h \approx 15.25$ meters. Therefore, the building should have a maximum height of approximately 15.25 meters. Considering that in civil construction each floor has an average height of 3 meters, this building corresponds to a building of approximately 5 floors.

Returning to the construction of the regular pentagon presented here, we see that point $C$ in Figure 5 is one of the intersection points of the line containing points $B$ and $O$ with the circumference with center $B$ and radius equal to $AB$. Considering then the other intersection point, say $C'$, we obtain the so-called golden rectangle determined by the consecutive vertices $A$, $O$ and $C'$, with sides measuring 1 and $(1 + \sqrt{5})/2$ (see [1] and [4]). An important highlight of this rectangle regarding the construction of the regular pentagon is the construction of Luca Pacioli for the regular icosahedron, described as follows (see [3], p. 8). It is obtained from the intersection of three perpendicular golden rectangles, arranged in particular symmetry position as shown in Figure 8a. For each of the twelve vertices of these perpendicular golden rectangles there are five other vertices which, together with the initial vertex, form five equilateral triangles. Moreover, these five points are coplanar and determine a regular pentagon. The polyhedron determined by all twelve vertices is the regular icosahedron (Figure 8b), which is one of the unique five regular polyhedra discovered by ancient greek mathematicians.

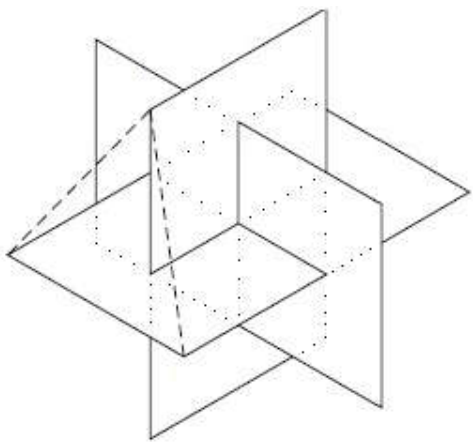 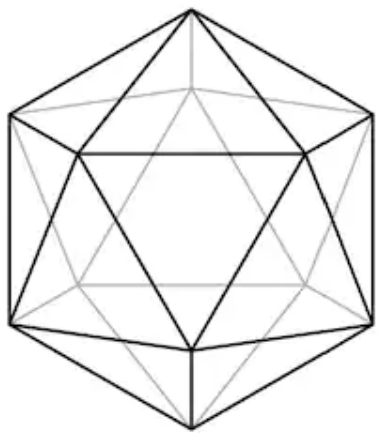

Figura 8a            Figura 8b

## Conclusion

In the case of the equilateral triangle (Figure 2), the mediatrix of segment $OA$ was taken. An analogous procedure can be done by taking the mediatrix of segment $OA'$. These two simultaneous procedures produce the regular hexagon $\wp_6$. Similarly, in the case of the regular pentagon (Figure 5), point $D$ was taken at segment $OA$ among the four possible intersections of the circumference with center $O$ and radius $OC$. The same procedure could be done using the intersection with $OA'$. These two simultaneous choices produce the decagon regular $\wp_{10}$, while all four simultaneous alternatives produce the regular icosagon $\wp_{20}$. Each of these cases leads us to the trigonometric ratios of angles 72°, 36° and 18°, as well as the angle 54°, which is the complementary angle of 36°. Moreover, as a stimulus to the reader, remembering that an angle inscribed on a circumference is half of its corresponding central angle (Figure 2), we even have the construction of the angle 9° and its complementary angle 81°, further enriching the application possibilities. Also, using the addition formulas, it is possible to calculate the trigonometric ratios of the angle 3° from the trigonometric ratios of the angles 18°, 30° and 45°.

A natural curiosity is to know what happens about the effective construction of $\wp_n$ if $n$ is greater than 5. We can separate the question into two situations, depending on the parity of $n$. If $n$ is even, say $n = 2k$, then the previously used method, considering the mediator of side $l_k$, can be applied assuming that the result for the regular polygon $\wp_k$ is already known. Otherwise, if $n$ is odd, then for $n = 7$ it is known that the method cannot be applied. In the general case, an answer was given by Gauss depending on the prime numbers factorization of $n$: $\wp_n$ is constructible if, and only, if $n = 2^r \cdot p_1 \cdots p_k$, where $r \in \mathbb{N}$ and $p_1, \cdots, p_k$ are distinct prime numbers of the form $p_i = 2^{2^{s_i}} + 1$, $1 \leq i \leq k$, $s_i \in \mathbb{N}$ (see for example [2]). Prime numbers of this form are known as Fermat prime numbers. Only the first five Fermat prime numbers are known, namely: 3, 5, 17, 257 and 65,537. It is known that several numbers of the form $2^{2^n} + 1$, $n \in \mathbb{N}$, are not primes. For example, Euler was the first person to warn about this, showing that the number $2^{2^5} + 1$ has 641 as prime factor and, therefore, is not prime.